\documentclass[12pt,reqno]{amsart}
\usepackage{amssymb}
\usepackage{cases}
\usepackage{bbm}
\usepackage{mathrsfs}
\usepackage{graphicx}
\usepackage{amsfonts}
\usepackage{enumerate}

\usepackage{amssymb, amsmath, amsfonts, latexsym}
\usepackage{enumerate}
\usepackage{ulem}
\usepackage{color}

\newtheorem{thm}{Theorem}[section]
\newtheorem{cor}[thm]{Corollary}
\newtheorem{lem}[thm]{Lemma}
\newtheorem{prop}[thm]{Proposition}
\newtheorem{example}[thm]{Example}
\newtheorem{remarks}[thm]{Remark}
\newtheorem{defn}[thm]{Definition}
\newtheorem{hyp}[thm]{Hypothesis}
\numberwithin{equation}{section}

\date{}

\newcommand{\ee}{\mathbb{E}}

\newcommand{\nn}{\mathbb{N}}
\newcommand{\rr}{\mathbb{R}}
\newcommand{\pp}{\mathbb{P}}

\def\BB{\mathcal B}

\def\FF{\mathcal F}

\def\MM{\mathcal M}

\def\<{\langle}
\def\>{\rangle}

\def\d"{^{\prime\prime}}

\def\bhyp{\begin{hyp}}
\def\nhyp{\end{hyp}}

\def\beq{\begin{equation}}
\def\nneq{\end{equation}}

\def\bdef{\begin{defn}}
\def\ndef{\end{defn}}

\def\bthm{\begin{thm}}
\def\nthm{\end{thm}}

\def\bprop{\begin{prop}}
\def\nprop{\end{prop}}

\def\brmk{\begin{remarks}}
\def\nrmk{\end{remarks}}

\def\bexa{\begin{example}}
\def\nexa{\end{example}}

\def\blem{\begin{lem}}
\def\nlem{\end{lem}}

\def\bcor{\begin{cor}}
\def\ncor{\end{cor}}

\def\bexe{\begin{exe}}
\def\nexe{\end{exe}}

\def\bprf{\begin{proof}}
\def\nprf{\end{proof}}

\def\dsp{\displaystyle}

\def\bdes{\begin{description}}
\def\ndes{\end{description}}

\title[Large deviations for mean-field Gibbs measures]
{Large deviations for empirical measures of mean-field Gibbs measures}
\author{Wei Liu}
\address{Wei Liu, School of Mathematics and Statistics, Wuhan University, Wuhan, Hubei 430072, PR China;
Computational Science Hubei Key Laboratory, Wuhan University, Wuhan, Hubei 430072, PR China.}
\thanks{The first author is supported by the CSC and NSFC(11731009, 11571262)}
\email{wliu.math@whu.edu.cn}

\author{Liming Wu}
\address{Liming Wu, Laboratoire de Math\'ematiques Blaise Pascal, CNRS-UMR
6620, Universit\'e Clermont-Auvergne, 3 Place Vasarely, 63178 Aubi\`ere, France.}
\email{Li-Ming.Wu@math.univ-bpclermont.fr}

\date{}

\begin{document}

\begin{abstract} In this paper, we show that   the empirical measure of mean-field model satisfies the large deviation principle with respect to the weak convergence topology or the stronger Wasserstein metric, under the strong exponential integrability condition on the negative part of the interaction potentials. In contrast to the known results we prove this without any continuity or boundedness condition on the interaction potentials. The proof relies mainly on the law of large numbers and the exponential decoupling inequality of de la Pe\~na for $U$-statistics.
\end{abstract}
\maketitle



\section{Introduction}
In this paper, we consider  a mean-field interacting system of $n$ particles at positions $x_1,\cdots,x_n$ in a separable and complete metric space $(S,\rho)$ (say Polish space) confined by a potential $V: S\rightarrow (-\infty, +\infty]$.
The interaction between the particles is given by a family of Borel-measurable interaction potentials $W^{(k)}: S^k\to (-\infty, +\infty]$ between $k$-particles,
where $2\le k\le N$, and $N\ge 2$ is fixed (necessarily $n\ge N$).
 The mean-field Hamiltonian
or energy functional $H_n : S^n\rightarrow (-\infty, +\infty]$ is given by
\begin{equation}\label{Hal} \aligned H_n(x_1,\cdots,x_n)&:=\sum^{n}_{i=1}V(x_i)+n\sum_{k=2}^N U_n(W^{(k)})\\
\endaligned
\end{equation}
where
\begin{equation}
 U_n(W^{(k)}) = \frac 1{|I_n^k|} \sum_{(i_1,\cdots,i_k)\in I_n^k}  W^{(k)}(x_{i_1},\cdots, x_{i_k})
\nneq
is the $U$-statistic of order $k$,
$$I_n^k:=\{(i_1,\cdots,i_k)\in \nn^k|i_1, \cdots, i_k \text{ are different }, 1\le i_1,\cdots, i_k\le n\}$$
and $|I_n^k|$ denotes the number of elements in $I_n^k$ (equal to $n !/(n-k) !$).

The mean-field Gibbs probability measure $P_n$ on $S^n$ is defined by
\begin{equation}
\label{Gibbs} \aligned dP_n(x_1,\cdots,x_n)&:=\frac{1}{Z_n}\exp(-H_n(x_1,\cdots,x_n))m(dx_1)\cdots m(dx_n)\\
\endaligned
\end{equation}
where $m$ is some nonnegative $\sigma$-finite measure on $S$ equipped with the Borel $\sigma$-field $\BB(S)$, and
\begin{equation}
\label{Zn} Z_n:=\int_{S}\cdots\int_{S}\exp(-H_n(x_1,\cdots,x_n))m(dx_1)\cdots m(dx_n)
\end{equation}
is the normalization constant (called {\it partition function}).

When $N=2$, this model is called mean-field of {\it pair interaction}, and when $N>2$, it is called mean-field of {\it many-bodies interaction}.

 The main objective of this paper is to study the large deviations of
the empirical measure
$$L_n(x^n;\cdot):=\frac1n\sum_{i=1}^{n}\delta_{x_i}(\cdot)$$
of configuration $x^n=(x_1,\cdots,x_n)$, under the mean-field measure $P_n$.  We will simply denote $L_n(x^n;\cdot)$ by $L_n$ when there is no likelihood of confusion.

In the case $S=\rr^d$, $P_n$ is just the equilibrium  state (or the invariant probability measure) of the system of $n$ interacting particles described by:
\begin{equation}
\label{SDE} dX^n(t)=\sqrt{2}d\mathrm{B}_t-\nabla H_n(X^n(t))dt,
\end{equation}
where $X^n(t):=(X_{1}^{n}(t),\cdots,X_{n}^{n}(t))^T$ ($\cdot^T$ means the {\it transposition}) takes values in $(\rr^{d})^n$, $\mathrm{B}_t=(B_t^1,\cdots,B_t^n)$, $B_t^1,\cdots,B_t^n$ are $n$ independent Brownian motions taking values in $\rr^d$. It is well-known that when $n$ goes to infinity,  $L_n(X^n(t); \cdot) $ converges to the solution of the nonlinear McKean-Vlasov equation (the so-called propagation of chaos), under quite general condition \cite{[Sz]}.

A classical problem is to establish conditions for the existence of a macroscopic
limit of the empirical measures $L_n$ as the number of particles $n\rightarrow+\infty.$ It is well-known that the large deviation principle (LDP in short) provides a strong exponential concentration with the speed $n$ in terms of some explicit rate function, which is very useful for the study of the macroscopic limit and microscopic phenomena in statistical mechanics. In the case of pair interaction (i.e. $N=2$),
L\'eonard established for the first time in \cite[1987]{[LC]} the LDP for the empirical measure $L_n$ under the Gibbs measure $P_n$ in the weighted weak convergence topology, when $\nu\to \iint W^{(2)}(x,y) d\nu(x)d\nu(y)$ is continuous in some appropriate topology and it is bounded by some weighted function satisfying the strong exponential integrability condition.
 By means of the weak convergence approach developed in Dupuis and Ellis \cite{[DE]}, Dupuis et al. established in \cite{[DLR]} an LDP in the case of pair interaction, by assuming that $W^{(2)}$ is lower bounded and lower semi-continuous (l.s.c. in short), which generalized the result obtained in \cite{[LC]}. For more results in this field the reader is referred to \cite{[BRJ],[BRJ1], [CGZ],[DG],[ZDG],[LS],[RS]} and the references therein.

The purpose of this paper is to establish the LDP for the empirical measures $L_n$  under a more general condition that $W^{(k)}$, $2\le k\le N$ are only measurable and their negative parts $W^{(k),-}$ satisfy the strong exponential integrability condition, which generalizes the previous results in \cite{[DLR]} and \cite{[LC]}. We first obtain the LDP with respect to the weak convergence topology, then with respect to the Wasserstein metric by using Sanov's theorem for the Wasserstein metric established by Wang et. al in \cite{[WWW]}. Our main tools are the law of large numbers (LLN in short) for the $U$-statistics and an exponential inequality for $U$-statistics issued of de la Pe\~na decoupling inequality.

The paper is organized as follows. In the next section, we will first briefly introduce some notations and definitions concerning the  LDP, and then present our main results. The proofs are presented in the third section.

\section{Main result}

\subsection{Preliminaries}\label{LDP}
We recall the definition of a rate function on a Polish space $S$ and the LDP for a sequence of probability measures on $(S,\mathcal{B}(S))$.
\bdef[Rate function]\label{rate function} $I$ is said to be a rate function on $S$ if it
is a lower semi-continuous function from $S$ to
$[0,\infty]$ (i.e., for all $L\ge0$, the level set $[I\le L]$ is
closed). $I$ is said to be a good rate function if it is inf-compact, i.e. $[I\le L]$ is
compact for any $L\in\rr$. \ndef
A consequence of a rate function being good is that its infimum is achieved over any non-empty closed set.

We denote by $\mathcal{M}_1(S)$ the space of probability measures on $S$.
\bdef[LDP]\label{LDP} Let $\{\nu_{n}\}_{n\in\nn}$ be a sequence of probability measures in $\mathcal{M}_1(S).$
\begin{enumerate}[$(a)$] \item
 $\{\nu_{n}\}_{n\in\nn}$ is said to satisfy the large deviation lower
bound with the speed $n$ and a rate function $I$ if for any open subset $G\in\mathcal
{B}(S)$, \beq\label{lower}
l(G):=\liminf_{n\to+\infty}\frac1n\log\nu_n(G)\ge-\inf_{\nu\in
G}I(\nu) ;\nneq

 \item $\{\nu_{n}\}_{n\in\nn}$ is said to satisfy the large deviation upper bound
with the speed $n$ and a rate function $I$ if for any closed subset $F\in\mathcal
{B}(S)$, \beq\label{upper}
u(F):=\limsup_{n\to+\infty}\frac1n\log\nu_n(F)\le-\inf_{\nu\in
F}I(\nu) ;\nneq

\item $\{\nu_{n}\}_{n\in\nn}$ is said to satisfy the large deviation principle
with the speed $n$ and a rate function $I$ if both (a) and (b) hold, and $I$ is good.
\end{enumerate}\ndef
The LDP characterizes the exponential concentration behavior, as $n\to+\infty$, of a sequence of probability measures $\{\nu_{n}\}_{n\in\nn}$ in terms
of a rate function. This characterization is via asymptotic upper and lower exponential bounds on the values that $\nu_{n}$
assigns to measurable subsets of $S$.

\subsection{Main results}
Throughout this paper, we assume that \beq\label{HV}C:=\int_{S}\exp(-V(x))m(dx)<+\infty.\nneq

Let \beq\label{V}\alpha(dx):=\frac1C\exp(-V(x))m(dx)\nneq be the probability measure on $S$,
then the mean-field Gibbs probability measure $P_n$ can be rewritten as
\begin{equation}
\label{Gibbs2} dP_n(x_1,\cdots,x_n)=\frac{1}{\widetilde{Z}_n}\exp\left(-n \sum_{k=2}^N U_n(W^{(k)})\right)\alpha^{\otimes n}(dx_1,\cdots, dx_n),
\end{equation}
where $\widetilde{Z}_n:=\frac{Z_n}{C^n}.$ Without interaction (i.e. $W^{(k)}=0$ for all $k$), $P_n=\alpha^{\otimes n}$, i.e. the $n$ particles are
free and identically distributed with law $\alpha$.

Given a probability measure $\mu\in\mathcal{M}_1(S)$, the relative entropy of $\nu$ with respect to $\mu$ is defined by
\beq \label{entropy}H(\nu|\mu)=\left\{
    \begin{array}{ll}
      \int_{S}\frac{d\nu}{d\mu}(x)\log\frac{d\nu}{d\mu}(x)\mu(dx), &\ \hbox{if $\nu \ll \mu$;} \\
      +\infty\ , &\ \hbox{otherwise.}
    \end{array}
  \right.
\nneq

For any probability measure $\nu\in\mathcal{M}_1(S)$ such that $W^{(k),-}:=(-W^{(k)})\vee0\in L^1(\nu^{\otimes k})$ for all $2\le k\le N$, we define
\beq \label{Jk}\mathcal{W}^{(k)}(\nu):=\int_{S^k}W^{(k)}(x_1,\cdots,x_k)d\nu^{\otimes k}(x_1,\cdots, x_k)\in(-\infty, +\infty].\nneq
The free energy of the state $\nu$ is given by
\beq \label{rfk}H_W(\nu):=\left\{
\begin{array}{ll}
      H(\nu|\alpha)+\dsp \sum_{k=2}^N\mathcal{W}^{(k)}(\nu), &\text{if }\  H(\nu|\alpha)<+\infty \text{ and }\\
&\quad\quad W^{(k),-}\in L^1(\nu^{\otimes k}), 2\le k\le N; \\
      +\infty\ , &\text{otherwise.}
    \end{array}
  \right.
\nneq
Without loss of generality we may and will assume that $W^{(k)}$ is symmetric,  i.e.
$$W^{(k)}(x_{\sigma(1)}, \cdots, x_{\sigma(k)})= W^{(k)}(x_1,\cdots,x_k)$$ for any $(x_1,\cdots,x_k)\in S^k$ and any permutation $\sigma$ on $\{1,\cdots,k\}$.
We make the following assumption on the interaction potentials $(W^{(k)})_{2\le k\le N}$:

\medskip
\textbf{(A1)} {\it For each $2\le k\le N$, the function $W^{(k)}: S^k\rightarrow(-\infty,+\infty]$ is symmetric, measurable; its positive part $W^{(k),+}$  satisfies
\beq\label{HWa}
H(\nu|\alpha) + \int_{S^k} W^{(k),+}(x_1,\cdots,x_k)d\nu^{\otimes k}(x_1,\cdots, x_k)<+\infty \ \text{ for {\bf some} } \nu\in\MM_1(S)
\nneq
and its negative part $W^{(k),-}$ satisfies the following {\bf strong exponential integrability condition}
\beq\label{HW}\mathbb{E}[\exp(\lambda W^{(k),-}(X_1,\cdots,X_k))]<+\infty, \ \forall \lambda>0\nneq
 where $X_1,\cdots,X_k$ are i.i.d. random variables of the common law $\alpha$ defined in (\ref{V}).}
\medskip

\begin{remarks}\label{A}{\rm  \begin{enumerate}[$(1)$]
\item The simplest condition for (\ref{HWa}) is: there is some measurable subset $F$ of $S$ with $\alpha(F)>0$ such that $1_{F^k}W^{(k),+}$ is  $\alpha^{\otimes k}$-integrable. In fact one can take $\nu=h \alpha$, where the density $h: S\to\rr^+$ is bounded, with support contained in $F$.

\item Under the exponential integrability condition (\ref{HW}), if $H(\nu|\alpha)<+\infty$, then $W^{(k), -}\in L^1(\nu^{\otimes k})$. In fact, for any $\lambda>0,$ by Donsker-Varadhan variational formula (see \cite[Lemma 1.4.3.(a)]{[DE]} in the bounded case) and Fatou's lemma (by approximating $W^{(k),-}$ with $W^{(k),-}\wedge L$, $L\uparrow +\infty$),
$$\aligned
\lambda&\int_{S^k}W^{(k), -}(x_1,\cdots,x_k)d\nu^{\otimes k}(x_1,\cdots, x_k)\\
&\le H(\nu^{\otimes k}|\alpha^{\otimes k})+\log\int_{S^k}e^{\lambda W^{(k), -}(x_1,\cdots,x_k)}d\alpha^{\otimes k}(x_1,\cdots, x_k)\\
&=k H(\nu|\alpha)+\log\int_{S^k}e^{\lambda W^{(k), -}(x_1,\cdots,x_k)}d\alpha^{\otimes k}(x_1,\cdots, x_k)<+\infty.\\
\endaligned$$

\item When $S=\rr^d$, $N=2$ and $\int_{\rr^d} e^{-V} dx<+\infty$,
 our assumption {\bf (A1)} is satisfied in the following two situations:

 \begin{enumerate}[(a)]
\item $W^{(2)}(x,y)= \frac b{|x-y|^\beta}$ with $\beta<d$, $b>0$ (Coulomb potential if $\beta=1$);

\item $W^{(2)}(x,y)= - b\log |x-y|$ with $b>0$ and $\int |x|^p e^{-V}dx<+\infty$ for all $p>0$.

(This  interaction potential appears in random matrices, and the result in \cite{[DLR]} does not apply for this example because of the lower boundedness assumption therein.)
 \end{enumerate}
In fact in both  cases, $W^{(2),-}$ satisfies the strong exponential integrability condition (\ref{HW}). Taking $F=\{x\in\rr^d; |x|\le R, V(x)\ge -L\}$ for $R,L>0$ large enough, we see that $1_{F^2}W^{(2),+}$ is $dxdy$-integrable, then $\alpha^{\otimes 2}$-integrable. So the condition (\ref{HWa}) is verified, by part (1) of this remark.
 \end{enumerate}
}
\nrmk

Now we present our first main result,  whose proof is given in the next section.

\bthm\label{main} Under assumption {\bf (A1)}, $H_W$ is inf-compact on $\mathcal{M}_1(S)$ equipped with the weak convergence topology, and
\beq\label{main1b}
-\infty<\inf_{\mu\in\mathcal{M}_1(S)}H_W(\mu)<+\infty.
\nneq
 Moreover the sequence of probability measures $\{P_n(L_n\in\cdot)\}_{n\ge N}$ satisfies the LDP on $\mathcal{M}_1(S)$ equipped with the weak convergence topology, with speed $n$ and the good rate function
 \beq\label{rfb} I_W(\nu):=H_W(\nu)-\inf_{\mu\in\mathcal{M}_1(S)}H_W(\mu), \ \nu\in\mathcal{M}_1(S).\nneq
Furthermore, for any $\nu$ such that $\nu\ll \alpha$ and $W^{(k),-}\in L^1(\nu^{\otimes k}),\ \forall 2\le k\le N$, the rate function
$I_W(\nu)$ can be identified as
 \beq\label{rfc}I_W(\nu)=\lim_{n\to+\infty}\frac1nH(\nu^{\otimes n}|P_n).\nneq
 \nthm

\brmk\label{rm1}{\rm Notice that we have removed  the lower semi-continuity and lower boundness conditions on  $W^{(k)}$ in \cite{[DLR]}. Our result generalizes the known results in L\'eonard \cite{[LC]}, Dupuis et al. \cite{[DLR]}.

}\nrmk

\brmk\label{rm2}{\rm To see the main difficulty in this LDP, let us proceed naively in the case of pair interaction: when $W(x,y):=W^{(2)}(x,y)$ is bounded and continuous, the $U$-statistic
$$
U_n(W)=\frac 1 {n(n-1)} \sum_{1\le i\ne j \le n} W(x_i,x_j)
$$
is very close to $\iint W(x,y) L_n(dx)L_n(dy)$ which is continuous in $L_n$ in the weak convergence topology (see the proof of Lemma \ref{LDP1} below).
So in that case the LDP follows from the LDP of $L_n$ under $\alpha^{\otimes \nn}$ (Sanov theorem) and Varadhan's Laplace lemma, as shown by L\'eonard \cite{[LC]}. When $W$ is bounded and only measurable,  we do not know whether the functional $\nu\to \iint W(x,y) \nu(dx)\nu(dy)$ is continuous in the (non-metrizable) $\tau$-topology, whereas the Sanov theorem still holds in the $\tau$-topology. The continuity of the last functional w.r.t. some appropriate topology is a basic assumption in  \cite{[LC]}.
}
\nrmk

\brmk\label{rm2}{\rm Since $H_W$ is inf-compact by Theorem \ref{main}, there is at least one minimizer. From the point of view of statistical physics, $H_W$ is an entropy or free energy associated to the nonlinear McKean-Vlasov equation. The uniqueness of the minimizer means that there is no phase transition for the mean-field. Below we recall some works on the uniqueness in the case of pair interaction.

For the uniqueness of the minimizer, it is sufficient to prove that $H_W$ is strictly convex along some path $(\nu_t)_{t\in[0,1]}$ connecting $\nu_0$ to $\nu_1$, for any two probability measures $\nu_0,\nu_1$. Let $\nu_*$ be a minimizer of $H_W$. Then the critical equation for the minimizer is
\beq\label{CE} \nu_*(dx)=\exp\left(-V(x)-   2\pi_{\nu_*} W^{(2)} (x)\right) m(dx)/C,
\nneq
where $$\pi_{\nu_*} W^{(2)} (x):=\int_{S} W^{(2)}(x,y) \nu_*(dy),$$
and $$C:=\int_{S}\exp\left(-V(x)-2\pi_{\nu_*}W^{(2)}(x)\right)m(dx)$$ is the normalization constant.

If $S=\rr^d$, the critical equation above is equivalent to the following stationary equation of the nonlinear McKean-Vlasov equation:
\beq\label{SE} \triangle\nu_*+\nabla\cdot(\nu_*\nabla V)+\nabla\cdot[(2\nabla W^{(2)}\ast\nu_*)\cdot\nu_*]=0,
\nneq
where the symbols $\nabla$ and $\nabla\cdot$ denote the gradient operator and divergence operator respectively. For the uniqueness of the solution of (\ref{SE}), the reader is referred to McCann \cite{[MC]} and Carrilo et al. \cite{[CMV]}. These authors showed that $H_W$ is strictly displacement convex (i.e. along the $W_2$-geodesic) under various sufficient conditions on the convexity of the confinement potential $V$ and the pair interaction potential $W^{(2)}$.

}
\nrmk

We also consider $\mathcal{M}_1(S)$ equipped with the Wasserstein topology, which is much stronger than the weak convergence topology. The {\bf $L^p$-Wasserstein distance} ($p\ge1$) with respect to the metric $\rho$ on $S$, between any two probability measures $\mu$ and $\nu$ on $S$, is defined by
\beq\label{Wa} W_p(\mu,\nu)=\inf_{\xi\in\Pi(\mu,\nu)}\left(\iint_{S\times S}\rho^p(x,y)\xi(dx,dy)\right)^{1/p},
\nneq
where $\Pi(\mu,\nu)$ is the set of all probability measures on $S\times S$ with marginal distribution $\mu$ and $\nu$ respectively (say couplings of $\mu$ and $\nu$).

The {\bf Wasserstein space of order $p$} is defined as $$\MM_1^p(S)=\left\{\mu\in\MM_1(S); \int_{S}\rho^p(x,x_0)\mu(dx)<+\infty\right\},$$
 where $x_0$ is some fixed point of $S$. It is known that $W_p$ is a finite distance on $\MM_1^p(S)$ and $(\MM_1^p(S), W_p)$ is a Polish space (see Villani \cite{[V1],[V2]}).

\bthm\label{mainW} Assume
\beq\label{mainW1}\int_{S}\exp\{\lambda \rho^p(x,x_0)\}\alpha(dx)<+\infty,\ \ \forall \lambda>0,\nneq
for some (hence for any) $x_0\in S$.
Under the assumption {\bf (A1)}, the sequence of probability measures $\{P_n(L_n\in \cdot)\}_{n\ge N}$ satisfies the LDP on $(\mathcal{M}_1^p(S),W_p)$ with speed $n$ and the good rate function $I_W$ defined in (\ref{rfb}).
\nthm

\brmk\label{rm3}{\rm In \cite{[DLR]}, Dupuis et al. imposed the following non-explicit condition for the LDP result above when $S=\rr^d$ and $N=2$: there exists a lower-semicontinuous function $\phi:\rr_{+}\to\rr$ with $$\lim_{s\to+\infty}\frac{\phi(s)}{s}=+\infty,$$
such that for every $\mu\in\MM_1(\rr^d)$,
$$\int_{\rr^d}\phi(|x|^p)\mu(dx)\le \inf_{\xi\in\Pi(\mu,\mu)}\left\{H(\xi|\alpha^{\otimes2})+\int_{\rr^d\times \rr^d}W(x,y)\xi(dxdy)\right\}.$$
}
\nrmk

\section{Proof of the main results}

Let $P_n^*$ be the measure by removing the normalization constant $\widetilde{Z}_n$ from $P_n$ presented in ($\ref{Gibbs2}$), i.e.,
\beq\label{PP} dP_n^*(x_1,\cdots,x_n):=\exp\left(-n\sum_{k=2}^N U_n(W^{(k)})\right)\alpha^{\otimes n}(dx_1, \cdots, dx_n).\nneq

\subsection{Large deviation (LD in short) lower bound for $P_n^*$}
First we present the law of large numbers of the $U$-statistic (see \cite[Corollary 3.1.1]{[KB]}). Let $X_1, X_2,\cdots$ be a sequence of  i.i.d. random variables in
a measurable space $(S, \mathcal{B}(S)).$ Let $\Phi:S^k \rightarrow \rr$ be a symmetric and measurable function of $k$ ($k\ge 2$) variables.
\blem\label{LLNU}\cite[Korolyuk and Borovskich]{[KB]} Assume that \beq\label{U3} \mathbb{E}|\Phi(X_1, \cdots, X_k)|<+\infty,
\nneq
then \beq\label{U4} U_n(\Phi)\rightarrow \mathbb{E}\Phi(X_1,\cdots, X_k)
\nneq
as $n\rightarrow +\infty$ with probability $1$.
\nlem
\bprf For the sake of completeness, we re-present the simple proof in \cite{[KB]}.

Let $\Pi_n$ be the set of all permutations of
$\{1,\cdots, n\}$ and
 $\mathfrak{B}_n$ the $\sigma$-algebra defined by
$$\mathfrak{B}_n:=\sigma\left\{B_n\times C_n|C_n\in\mathcal{B}(S^{[n+1,+\infty)}), B_n\in\mathcal{B}(S^{n}),\pi1_{B_n}=1_{B_n}, \forall\pi \in\Pi_n\right\}. $$
The $\sigma$-algebra $\mathfrak{B}_n$ remains unchanged under any permutation in $\Pi_n$, and $\mathfrak{B}_n\supseteq\mathfrak{B}_{n+1}$ for every $n\ge1.$

For any
$(i_1, \cdots, i_k)\in I_n^k,$ by (\ref{U3}) we have
$$\mathbb{E}[\Phi(X_{i_1}, \cdots, X_{i_k})|\mathfrak{B}_n]=\mathbb{E}[\Phi(X_1, \cdots, X_k)|\mathfrak{B}_n],$$ which yields $$U_n(\Phi)=\mathbb{E}[\Phi(X_1, \cdots, X_k)|\mathfrak{B}_n].$$
According to the limit theorem for reversed martingales and the $0$-$1$ law for $\mathcal{B}_{\infty}=\bigcap_{n\ge1}\mathcal{B}_n$,
$$U_n(\Phi)\xrightarrow{a.s.}\mathbb{E}[\Phi(X_1, \cdots, X_k)|\mathfrak{B}_{\infty}]=\mathbb{E}\Phi(X_1, \cdots, X_k).$$
\nprf

We have the following LD lower bound for the empirical measure $L_n$ under $P_n^*.$
\bprop\label{LLDP1} Without any integrability condition on $(W^{(k)})_{2\le k\le N}$,
the following large deviation lower bound holds for $\{P_n^*\{L_n\in\cdot\}\}_{n\ge N}$:
\beq\label{LLDP1a}l^{*}(G):=\liminf_{n\to+\infty}\frac1n\log P_n^*\{L_n\in G\}\ge-\inf\{H_W(\nu)|\nu\in G, W^{(k)}\in L^1(\nu^{\otimes k}),\ 2\le k\le N\},
\nneq
for
any open subset $G$ of $\mathcal{M}_1(S)$. In particular, we have
\beq\label{LZ}\liminf_{n\to+\infty}\frac1n\log\widetilde{Z}_n\ge -\inf\{H_W(\nu)|\nu\in \mathcal{M}_1(S), W^{(k)}\in L^1(\nu^{\otimes k}),\ 2\le k\le N\}.
\nneq
\nprop

\bprf Since (\ref{LZ}) is obtained just by taking $G$ as $\mathcal{M}_1(S)$ in (\ref{LLDP1a}), we only need to prove (\ref{LLDP1a}).
For (\ref{LLDP1a}), it is enough to show that for any $\nu\in G$ such that $H(\nu|\alpha)<+\infty$ and $W^{(k)}\in L^1(\nu^{\otimes k}),\ 2\le k\le N$,
$$l^{*}(G)\ge-H_W(\nu).$$

Let $\mathcal{N}(\nu, \delta)$ be the open ball centered at $\nu$ in $\mathcal{M}_1(S)$ with radius $\delta$ in the L\'evy-Prokhorov metric $d_w$ such that $\mathcal{N}(\nu, \delta)\subset G$.
Introduce the events
$$
\aligned A_n&:= \{x^n=(x_1,\cdots, x_n)\in S^n |L_n=L_n(x^n,\cdot)\in \mathcal{N}(\nu, \delta)\},\\
B_n&:=\{(x_1,\cdots, x_n)|\frac{1}{n}\sum_{i=1}^n\log\frac{d\nu}{d\alpha}(x_i)\le H(\nu|\alpha)+\varepsilon\},\\
C_n&:=\{(x_1,\cdots, x_n)|\sum_{k=2}^N U_n(W^{(k)})\le \sum_{k=2}^N \mathcal{W}^{(k)}(\nu)+\varepsilon\}.
\endaligned$$
Then we have for any $\varepsilon>0,$
\beq\label{lld1}  \aligned P_n^*\{L_n\in \mathcal{N}(\nu, \delta)\}&\ge \int_{A_n}\left(\frac{d\nu^{\otimes n}}{dP_n^*}(x_1,\cdots,x_n)\right)^{-1}d\nu^{\otimes n}(x_1, \cdots, x_n)\\
&=\int_{A_n}\exp\left(-\sum_{i=1}^{n}\log\frac{d\nu}{d\alpha}(x_i)\right)\exp\left(-n\sum_{k=2}^N U_n(W^{(k)})\right)\nu^{\otimes n}(dx_1, \cdots, dx_n)\\
&\ge \nu^{\otimes n}(A_n\cap B_n \cap C_n)\exp\left(-n[H(\nu|\alpha)+\varepsilon]-n[\sum_{k=2}^N \mathcal{W}^{(k)}(\nu)+\varepsilon]\right)      \\
&=\nu^{\otimes n}(A_n\cap B_n \cap C_n)\exp\left(-nH_W(\nu)-2n\varepsilon\right),  \\
\endaligned
\nneq

We claim that $\lim_{n\rightarrow +\infty}\nu^{\otimes n}(A_n\cap B_n \cap C_n)=1.$
Indeed, by the LLN, it is obvious that $$ \nu^{\otimes n}(A_n)\rightarrow1 \ \text{and}\ \nu^{\otimes n}(B_n)\rightarrow1 \ \text{as}\ n\rightarrow+\infty.$$
By the LLN of $U$-statistics in Lemma \ref{LLNU}, we also have $$\lim_{n\rightarrow +\infty}\nu^{\otimes n}(C_n)=1.$$

With this claim in hand, we immediately get from (\ref{lld1}) that
\beq\label{lld2} l^{*}(G)\ge\liminf_{n\rightarrow+\infty}\frac1n\log P_n^*\{L_n\in \mathcal{N}(\nu, \delta)\} \ge -H_W(\nu)-2\varepsilon,
\nneq
which completes the proof since $\varepsilon>0$ is arbitrary.
\nprf

\subsection{Decoupling inequality of de la Pe\~na and the key lemma}

We first recall  the decoupling inequality of de la Pe\~na \cite[1992]{[DELA]}.

\bprop\label{Decoupling2}\cite[de la Pe\~na]{[DELA]} Let $\{X_i\}_{i\ge1}$ be a family of i.i.d. random variables in a measurable space $(S, \mathcal{B}(S))$ and suppose that $(X_1^j, \cdots, X_n^j)_{j=1}^k$ are independent copies of $(X_1,\cdots,X_n)$. Let $\Psi$ be any convex increasing function on $[0,+\infty)$. Let $\Phi:S^k\rightarrow \rr$ be a symmetric and measurable function of $k$ variables such that
\beq\label{Dec12}\mathbb{E}|\Phi(X_1, \cdots, X_k)|<+\infty,
\nneq
then
\beq\label{Dec22}\mathbb{E}\Psi\left(\left| \sum_{(i_1,i_2,\cdots,i_k)\in I_n^k}\Phi(X_{i_1}, \cdots, X_{i_k})\right|\right)\le \mathbb{E}\Psi\left(C_k\left|\sum_{(i_1,i_2,\cdots,i_k)\in I_n^k}\Phi(X_{i_1}^1, \cdots, X_{i_k}^k)\right|\right),
\nneq where $C_2=8$ and $\dsp C_k=2^k \prod_{j=2}^k (j^j-1)$ for $k>2$.
 \nprop

Besides the decoupling inequality of de la Pe\~na above, we also require the following inequality.

\blem\label{U-Lap3} Let $1\le k\le n$, and $\{X_i^j;1\le i\le n, 1\le j\le k\}$ be independent random variables. For any $(i_1,\cdots,i_k)\in I_n^k$, let $\Phi_{i_1,\cdots,i_k}:S^k\rightarrow \rr$ be a measurable function of $k$ variables, then
\beq\label{Lap31}\aligned &\log\ee\exp\left(\frac{(n-k)!}{n!}\sum_{(i_1,\cdots,i_k)\in I_n^k}\Phi_{i_1,\cdots,i_k}(X_{i_1}^1,\cdots,X_{i_k}^k)\right)\\
&\le \frac{(n-k+1)!}{n!}\sum_{(i_1,\cdots,i_k)\in I_n^k}\log\ee\exp\left(\frac{1}{n-k+1}\Phi_{i_1,\cdots,i_k}(X_{i_1}^1,\cdots,X_{i_k}^k)\right).\\
\endaligned
\nneq
\nlem

\bprf For $k=1$, (\ref{Lap31}) is obviously  an equality. Next we prove this lemma by induction.
Assume that (\ref{Lap31}) is valid for $k-1$. Denote the left-hand side of (\ref{Lap31}) by $B_k$ and write $\sum_{(i_1,\cdots,i_{k-1})\in I_n^{k-1}}$ as $\sum_{I_n^{k-1}}$ for simplicity. We have
\beq\label{lap32}\aligned B_k&=\log\ee^{X^k}\{\ee[\exp(\frac{(n-k+1)!}{n!}\sum_{I_n^{k-1}}\sum_{i_k\notin \{i_1,\cdots, i_{k-1}\}}\frac1{n-k+1}\Phi_{i_1,\cdots,i_k}(X_{i_1}^1,\cdots,X_{i_k}^k))|X^k]\}.\\
\endaligned
\nneq

Given $X^k=(X^k_1,\cdots,X^k_n)$, let
\beq\label{Phi}\tilde{\Phi}_{i_1,\cdots,i_{k-1}}:=\frac1{n-k+1}\sum_{i_k:i_k\notin \{i_1,\cdots, i_{k-1}\}}\Phi_{i_1,\cdots,i_k}(X_{i_1}^1,\cdots, X_{i_{k-1}}^{k-1} ,X_{i_k}^k).\nneq
By the assumption of $(k-1)^{th}$ step, we get

\beq\label{lap33}\aligned B_k&\le \log\ee^{X^k}\left\{\exp\left(\frac{(n-k+2)!}{n!}\sum_{I_n^{k-1}}\log\ee[
\exp(\frac1{n-k+2}\tilde{\Phi}_{i_1,\cdots,i_{k-1}})|X^k]\right)\right\}\\
&=\log\ee^{X^k}\left\{\exp\left(\frac{(n-k+1)!}{n!}\sum_{I_n^{k-1}}\log\{\ee[
\exp(\frac1{n-k+2}\tilde{\Phi}_{i_1,\cdots,i_{k-1}})|X^k]\}^{n-k+2}\right)\right\}\\
&\le\frac{(n-k+1)!}{n!}\sum_{I_n^{k-1}}\log\ee^{X^k}\left\{\left[\ee\exp\left(\frac1{n-k+2}
\tilde{\Phi}_{i_1,\cdots,i_{k-1}}\right)|X^k\right]^{n-k+2}\right\}\\
\endaligned
\nneq
where the last inequality follows by the convexity of $X\to\log\ee e^X$ (a consequence of H\"older's inequality).

Next we deal with the logarithmic term in the last inequality above. Given $(i_1,\cdots, i_{k-1})$,
\beq\label{lap34}\aligned
&\ee^{X^k}\{[\ee\exp(\frac1{n-k+2}
\tilde{\Phi}_{i_1,\cdots,i_{k-1}})|X^k]^{n-k+2}\}\\
&=\ee^{X^k}\{[\ee\exp(\frac1{(n-k+2)(n-k+1)}
\sum_{i_k\notin \{i_1,\cdots, i_{k-1}\}}\Phi_{i_1,\cdots,i_k}(X_{i_1}^1,\cdots,X_{i_{k-1}}^{k-1},X_{i_k}^k))|X^k]^{n-k+2}\}\\
&\le \ee^{X^k}\left\{\left[\Pi_{i_k:i_k\notin \{i_1,\cdots, i_{k-1}\}}\ee[\exp(\frac1{n-k+2}
\Phi_{i_1,\cdots,i_k}(X_{i_1}^1,\cdots,X_{i_{k-1}}^{k-1},X_{i_k}^k))|X^k]\right]^{\frac{n-k+2}{n-k+1}}\right\}\\
&\ \ \ \ \ \ \ \ \ \ \ \ \ \ \ \ \ \ \ \ \ \ \ \ \ \ \ \ \ \ \ \ \  \ \ \ \ \ \ \ \ \ \ \ \ \ \ \ \ \ \ \ \ \ \ \ \ \ \ \ \ (\text{by H\"older's inequality})\\
&\le \ee^{X^k}\left\{\Pi_{i_k:i_k\notin \{i_1,\cdots, i_{k-1}\}}\ee[\exp(\frac1{n-k+1}
\Phi_{i_1,\cdots,i_k}(X_{i_1}^1,\cdots,X_{i_{k-1}}^{k-1},X_{i_k}^k))|X^k]\right\}\\
&\ \ \ \ \ \ \ \ \ \ \ \ \ \ \ \ \ \ \ \ \ \ \ \ \ \ \ \ \ \ \ \ \  \ \ \ \ \ \ \ \ \ \ \ \ \ \ \ \ \ \ \ \ \ \ \ \ \ \ \ \ (\text{by Jensen's inequality})\\
&=\Pi_{i_k:i_k\notin \{i_1,\cdots, i_{k-1}\}}\ee\exp(\frac1{n-k+1}
\Phi_{i_1,\cdots,i_k}(X_{i_1}^1,\cdots,X_{i_{k-1}}^{k-1},X_{i_k}^k))\\
&\ \ \ \ \ \ \ \ \ \ \ \ \ \ \ \ \ \ \ \ \ \ \ \ \ \ \ \ \ \ \ \ \  \ \ \ \ \ \ \ \ \ \ \ \ \ \ \ \ \ \ \ \ \ \ \ \ \ \ \ \ (\text{by the independence of }X_1^k, \cdots, X_n^k)\\
\endaligned
\nneq
Plugging (\ref{lap34}) into (\ref{lap33}), we get the desired inequality (\ref{Lap31}).

\nprf

Let $\{X_i\}_{i\ge1}$ be a family of i.i.d. random variables of law $\alpha$. Denote by $\Lambda_{n}(\cdot\ ;W^{(k)})$ the logarithmic moment generating function associated with the $U$-statistic of order $k$, i.e., for any $n\ge k\ge2$ and $\lambda>0,$ \beq\label{Lap2}\Lambda_{n}(\lambda;W^{(k)}):=\frac1n\log\mathbb{E}[\exp(\lambda nU_{n}(W^{(k)}))].
\nneq
Now we present our key lemma, which is a corollary of  the decoupling inequality of de la Pe\~na and Lemma \ref{U-Lap3}.
\blem\label{U-Lap} If $\mathbb{E}|W^{(k)}(X_1,\cdots,X_k)|<+\infty,$ then for any $2\le k\le n$ and $\lambda>0,$
\beq\label{Lap21}\Lambda_{n}(\lambda;W^{(k)})\le \frac1k\log\mathbb{E}[\exp(kC_k\lambda |W^{(k)}(X_1,\cdots,X_k)|)],
 \nneq where $C_k$ is defined as in Proposition \ref{Decoupling2}.
\nlem

\bprf Let $(X_1^j, \cdots, X_n^j)_{j=1}^k$ be independent copies of $(X_1,\cdots,X_n)$. By Lemma \ref{Decoupling2} (the decoupling inequality of de la P\~ena) and Lemma \ref{U-Lap3}, taking $\Phi_{i_1,\cdots,i_k}\equiv W^{(k)}$ for any $ (i_1,\cdots,i_k)\in I_n^k$, we get for any $\lambda>0$,
\beq\label{Lap35}\aligned &\Lambda_{n}(\lambda;W^{(k)})=\frac1n\log\mathbb{E}[\exp(\lambda n\frac{(n-k)!}{n!}\sum_{(i_1,\cdots,i_k)\in I_n^k}W^{(k)}(X_{i_1},\cdots,X_{i_k}))]\\
&\le\frac1n\log\mathbb{E}[\exp(\frac{(n-k)!}{n!}\sum_{(i_1,\cdots,i_k)\in I_n^k}\lambda nC_k |W^{(k)}|(X_{i_1}^1,\cdots,X_{i_k}^k)]\ \ \ \ (\text{by (\ref{Dec22})})\\
&\le \frac 1n \frac{(n-k+1)!}{n!} \sum_{(i_1,\cdots,i_k)\in I_n^k} \log \ee \exp(\frac{\lambda nC_k}{n-k+1}|W^{(k)}|(X_{i_1}^1,\cdots,X_{i_k}^k))\ \ \ \ (\text{by (\ref{Lap31})})\\
&= \frac{n-k+1}{n} \log\mathbb{E}\exp(\frac{\lambda nC_k}{n-k+1}|W^{(k)}|(X_{1},\cdots,X_{k}))\\
&\le  \frac1k\log\ee\exp(kC_k\lambda|W^{(k)}|(X_{1},\cdots,X_{k})),\\
\endaligned
\nneq
where the last inequality follows by the non-decreasingness of $a\to \frac 1a \log \ee e^{aX}$ on $(0,+\infty)$ and  the fact that $\frac{n}{n-k+1}\le k$ for all $n\ge k$.
\nprf

We have the following exponential approximation of the $U$-statistics.

\blem\label{U-exp}Assume that for any $\lambda>0, $
\beq\label{W2}\mathbb{E}[\exp(\lambda |W^{(k)}|(X_1,\cdots,X_k))]<+\infty.\nneq
Then there exists a sequence of bounded continuous functions $\{W_m^{(k)}\}_{m\ge1}$ such that for any $\delta>0,$
\beq\label{U-exp2}\lim_{m\rightarrow +\infty}\limsup_{n\rightarrow +\infty}\frac1n\log\mathbb{P}\{|U_n(W^{(k)})-U_n(W_m^{(k)})|>\delta\}=-\infty.
\nneq
\nlem
\bprf
For any function $W^{(k)}$ satisfying (\ref{W2}) and integer $m\in\nn^*$, considering the truncation $W^{(k),L}:=(-L)\vee (W^{(k)}\wedge L)$, we have by dominated convergence that
$$\log\mathbb{E}[\exp(m |W^{(k)}-W^{(k),L}|(X_1,\cdots,X_k))]\rightarrow 0$$ as $L\to+\infty$. Then we can choose $L=L(m)$ so that
 $$\log\mathbb{E}[\exp(m |W^{(k)}-W^{(k),L(m)}|(X_1,\cdots,X_k))]\le \frac 1m.
 $$ For $m\in\nn^*$ and $L=L(m)>0$ fixed, we can find a sequence of  bounded continuous functions $\{W^{(k),L}_l\}_{ l\ge1}$ on $S^k$ such that $W^{(k),L}_l(X_1,\cdots,X_k)\to W^{(k),L}(X_1,\cdots,X_k)$ in $L^1$  as $l$ goes to infinity, and $|W^{(k),L}_l(X_1,\cdots,X_k)|\le L$ (otherwise considering the truncation $(-L)\vee (W^{(k),L}_l\wedge L)$) for all $l\ge1$. Since
 $$\aligned\exp&(m [ |W^{(k)}-W^{(k),L}|(X_1,\cdots,X_k)+ |W^{(k),L}-W^{(k),L}_l|(X_1,\cdots,X_k)])\\
 &\le \exp(m  |W^{(k)}-W^{(k),L}|(X_1,\cdots,X_k)+2mL)\endaligned$$
for any $l\ge 1$, we have by the dominated convergence
 $$\aligned
 &\mathbb{E}[\exp(m [ |W^{(k)}-W^{(k),L}|(X_1,\cdots,X_k)+ |W^{(k),L}-W^{(k),L}_l|(X_1,\cdots,X_k)])]\\
 &\to
 \mathbb{E}[\exp(m  |W^{(k)}-W^{(k),L}|(X_1,\cdots,X_k))]
 \endaligned$$
as $l\to\infty$. Thus for $L=L(m)$, we can find $l=l(m)$ such that
$$
 \log\mathbb{E}[\exp(m [ |W^{(k)}-W^{(k),L}|(X_1,\cdots,X_k)+ |W^{(k),L}-W^{(k),L}_l|(X_1,\cdots,X_k)])]\le \frac 2m.$$
Setting $W^{(k)}_m=W^{(k),L(m)}_{l(m)}$ which is bounded and continuous, we have by the triangular inequality
$$
\log\mathbb{E}[\exp(m |W^{(k)}-W^{(k)}_m|(X_1,\cdots,X_k))]\le \frac 2m.
$$
For any $\lambda>0$, since
$$\log\mathbb{E}[\exp(\lambda |W^{(k)}-W_m^{(k)}|(X_1,\cdots,X_k))]\le\frac{\lambda}{m} \log\mathbb{E}[\exp(m|W^{(k)}-W_m^{(k)}|(X_1,\cdots,X_k))]$$
 for $m\ge \lambda$ by Jensen's inequality, we see that for the sequence of bounded and continuous functions $\{W^{(k)}_m\}_{m\ge1}$,
\beq\label{U-exp4}\log\mathbb{E}[\exp(\lambda |W^{(k)}-W_m^{(k)}|(X_1,\cdots,X_k))]\rightarrow 0,  \nneq
as $m\rightarrow+\infty$.

For any $\delta,\lambda>0,$ by Chebyshev's inequality we have
\beq\label{U-exp3}  \aligned \mathbb{P}\{|U_n(W^{(k)})-U_n(W_m^{(k)})|>\delta\}&\le e^{-n\lambda\delta}\mathbb{E}[\exp(\lambda nU_n(|W^{(k)}-W_m^{(k)}|))].\\
\endaligned
\nneq
Applying Lemma \ref{U-Lap}, we get
\beq\label{U-exp5}\frac1n\log\mathbb{P}\{|U_n(W^{(k)})-U_n(W_m^{(k)})|>\delta\}\le -\lambda\delta+\frac1k\log\mathbb{E}[\exp(kC_k\lambda |W^{(k)}-W_m^{(k)}|(X_1,\cdots,X_k))].
\nneq
Letting $m\rightarrow+\infty$ and using (\ref{U-exp4}), we get  the desired result (\ref{U-exp2}), since $\lambda>0$ is arbitrary.
\nprf

\subsection{LDP of $U$-statistics}
We begin with
\blem\label{LDP1} Let $(X_n)_{n\ge1}$ be a sequence of i.i.d. random variables valued in $S$, of common law $\alpha$.
Assume that $(W^{(k)})_{2\le k\le N}$ are measurable and satisfy the strong exponential integrability condition
\beq \label{LDP11}
\ee \left[\exp\left(\lambda |W^{(k)}(X_1,\cdots,X_k)|\right)\right]<+\infty, \ \forall \lambda>0,
\nneq
then $\{\pp((L_n, U_n(W^{(2)}),\cdots, U_n(W^{(N)}))\in \cdot)\}_{n\ge N}$ satisfies the LDP  on the product space $\mathcal{M}_1(S)\times\rr^{N-1}$, with good rate function $I$ defined by
\beq \label{rf2}I(\nu,z_2,\cdots, z_N):=\left\{
    \begin{array}{ll}
      H(\nu|\alpha), & \hbox{if $z_k=\mathcal{W}^{(k)}(\nu)$ for all $2\le k\le N$} \\
      +\infty\ , & \hbox{otherwise}
    \end{array}
  \right.
\nneq
for any $(\nu,z_2,\cdots, z_N)\in \mathcal{M}_1(S)\times \rr^{N-1}$.
\nlem

\bprf  Let $\{W^{(k)}_m\}_{m\ge1}, 2\le k\le N$ be the sequences of bounded continuous functions as in Lemma \ref{U-exp} such that for any $\lambda>0$,
 \beq\label{LDP111} \varepsilon(\lambda,m,k):=\log\int_{S^k} \exp(\lambda|W^{(k)}_m-W^{(k)}|)d\alpha^{\otimes k}\rightarrow 0,\ \text{as} \ m\rightarrow+\infty.\nneq

 For any $m\ge1$, let $$ \aligned f_m(\nu)&:=\left(\nu, \int_{S^2}W^{(2)}_md\nu^{\otimes 2},\cdots, \int_{S^N}W^{(N)}_md\nu^{\otimes N}\right),\\
  f(\nu)&:=\left(\nu, \int_{S^2}W^{(2)}d\nu^{\otimes 2},\cdots, \int_{S^N}W^{(N)}d\nu^{\otimes N}\right)\\
  \endaligned $$
and consider the following metric on the product space $\mathcal{M}_1(S)\times\rr^{N-1}$:
$$d((\nu_1,z_2,\cdots,z_N), (\tilde{\nu}_1,\tilde{z}_2,\cdots,\tilde{z}_N)):=d_w(\nu_1,\tilde{\nu}_1) +\sum_{k=2}^N|z_k-\tilde{z}_k|.$$
Note that $d(f_m(\nu),f(\nu))=\sum_{k=2}^N\left|\int_{S^k}(W^{(k)}_m-W^{(k)})d\nu^{\otimes k}\right|$. Below we separate the proof of the LDP into
three points.

{\bf 1) The continuity of $f_m(\nu)$.}
We first prove the continuity of $f_m(\nu)$ or equivalently that of $\nu\to \int W_m^{(k)} d\nu^{\otimes k}$ on $\MM_1(S)$ in the weak convergence topology£¬ for each $k=2,\cdots, N$.
 In fact let $\nu_n\to\nu$ in $(\MM_1(S), d_w)$. By Skorokhod's representation theorem (see \cite[Theorem 6.7]{[Bill]}), one can construct a sequence of $S$-valued random variables  $Y_n$ of law $\nu_n$, converging a.s. to $Y$ of law $\nu$. Let $(Y_n^{(i)}, n\ge 0; Y^{(i)})_{1\le i\le k}$ be independent copies of $(Y_n, n\ge 0; Y)$. Then $Y_n^{(i)}\to Y^{(i)}$, $a.s.$ too for each $i$, in other words,  $(Y_n^{(1)},\cdots,Y_n^{(k)})\to (Y^{(1)},\cdots,Y^{(k)}), a.s.$. Thus $\nu_n^{\otimes k} \to \nu^{\otimes k}$ weakly on $S^k$, which implies the continuity of the functional above.

{\bf 2) Exponentially good approximation  of $(L_n, U_n(W^{(2)}),\cdots, U_n(W^{(N)}))$ by $f_m(L_n)$.}

By (\ref{U-exp2})  in Lemma \ref{U-exp}, for any $\delta>0$,
$$\aligned
\limsup_{n\to\infty} \frac 1n \log &\pp\left(d((L_n, U_n(W_m^{(2)}),\cdots, U_n(W_m^{(N)})), (L_n, U_n(W^{(2)}),\cdots, U_n(W^{(N)})))>\delta\right)\\
&\to -\infty
\endaligned$$
as $m\to+\infty$, i.e. $(L_n, U_n(W_m^{(2)}),\cdots, U_n(W_m^{(N)}))$ is  an exponentially good approximation of $(L_n, U_n(W^{(2)}),\cdots, U_n(W^{(N)}))$.

Moreover $(L_n, U_n(W_m^{(2)}),\cdots, U_n(W_m^{(N)}))$  and $f_m(L_n)$ are exponentially equivalent, because the following uniform estimate holds:
$$
\aligned
|U_n(W^{(k)}_m)-\int W^{(k)}_m dL_n^{\otimes k}|
&\le |U_n(W^{(k)}_m) - \frac{|I_n^k|}{n^k}U_n(W^{(k)}_m)| + (1-\frac{|I_n^k|}{n^k})\|W^{(k)}_m\|_\infty\\
&\le 2\left(1-\frac{|I_n^k|}{n^k}\right)\|W^{(k)}_m\|_\infty\to 0\\
\endaligned
$$
as $n\to\infty$. Then as $m\to\infty$, $f_m(L_n)$ is an exponentially good approximation of  $(L_n, U_n(W^{(2)}),\cdots, U_n(W^{(N)}))$.

{\bf 3).}
 Therefore by Sanov's theorem (the LDP of $L_n$) and the theorem of approximation of LDP (\cite[Theorem 4.2.23]{[DZ]}), for the desired LDP it remains to show that  for any $L>0$,
\beq\label{LDP113}\sup_{\nu:H(\nu|\alpha)\le L} d(f_m(\nu), f(\nu))\to 0,\ \text{as} \ m\rightarrow+\infty.\nneq

In fact for any $\lambda>0,\ L>0$ and $\nu$ with $H(\nu|\alpha)\le L$, we have by Donsker-Varadhan variational formula (see \cite[Lemma 1.4.3.(a)]{[DE]}) and Fatou's lemma,
\beq\label{LDP112}\aligned
\int_{S^k} &|W^{(k)}_m-W^{(k)}|d\nu^{\otimes k}\\
&\le \frac{1}{\lambda}\left(H(\nu^{\otimes k}|\alpha^{\otimes k})+\log\int_{S^k}\exp(\lambda|W^{(k)}_m-W^{(k)}|)d\alpha^{\otimes k}\right)\\
&\le \frac{1}{\lambda}[kL+\varepsilon(\lambda,m,k)],\ \ 2\le k\le N.\\
\endaligned\nneq
where (\ref{LDP113}) follows for $\lambda>0$ is arbitrary and $\lim_{m\to\infty} \varepsilon(\lambda,m,k)=0$ for any $\lambda>0$.

\nprf

\subsection{Proof of Theorem \ref{main}}

\bprf[Proof of Theorem \ref{main}] We divide its proof into three steps.

{\bf Step 1. $W^{(k)}$ is upper bounded.}
 In this upper bounded case, $\ee e^{\lambda |W^{(k)}|(X_1,\cdots, X_k)}<+\infty$ for all $\lambda>0$  by $\textbf{(A1)}$. By Lemma \ref{LDP1}, under $\pp=\alpha^{\otimes \nn^*}$, $(L_n, U_n(W^{(2)}), \cdots, U_n(W^{(N)}))$ satisfies the LDP on $\MM_1(S)\times \rr^{N-1}$ with the rate function $I(\nu,z_2,\cdots, z_N)$ given by (\ref{rf2}). Since $\sum_{k=2}^N U_n(W^{(k)})$ is continuous in $(L_n, U_n(W^{(2)}), \cdots, U_n(W^{(N)}))$,  and
$$
\limsup_{n\to\infty}\frac 1n \log \ee \exp\left(-n p \sum_{k=2}^N U_n(W^{(k)})\right) <+\infty
$$
for some $p>1$ (in fact for all $p>1$) by Lemma \ref{U-Lap} and $\textbf{(A1)}$, we can apply the tilted LDP (Ellis \cite[Theorem II.7.2.]{Ellis}) to conclude that  $P_n((L_n, U_n(W^{(2)}), \cdots, U_n(W^{(N)}))\in\cdot)$ satisfies the LDP, with the rate function $\tilde I$ given by
$$\aligned
\tilde I&(\nu,z_2,\cdots, z_N) =I(\nu,z_2,\cdots, z_N)+ \sum_{k=2}^N z_k  -  \inf_{(\nu,z_2,\cdots, z_N)} [I(\nu,z_2,\cdots, z_N)+\sum_{k=2}^N z_k ]\\
&= \begin{cases} &H(\nu|\alpha) +\sum_{k=2}^N {\mathcal W}^{(k)}(\nu) - \inf\limits_{\nu'\in \MM_1(S), H(\nu'|\alpha)<+\infty}\{H(\nu'|\alpha) +\sum_{k=2}^N {\mathcal W}^{(k)}(\nu')\}\\
&\quad\quad\quad\quad \ \text{ if }\ H(\nu|\alpha)<+\infty,  z_k={\mathcal W}^{(k)}(\nu), 2\le k\le N\\
&+\infty,\quad \text{ otherwise.}
\end{cases}\\
&= \begin{cases} H_W(\nu) - \inf_{\MM_1(S)}H_W\ &\text{ if } \ H(\nu|\alpha)<+\infty,  z_k={\mathcal W}^{(k)}(\nu), 2\le k\le N\\
+\infty &\text{ otherwise.}
\end{cases}
\endaligned$$
Hence $P_n(L_n\in\cdot)$ satisfies the LDP with the rate function $I_W$, by the contraction principle. Notice that $H_W = I_W +\inf_{\MM_1(S)}H_W$ is inf-compact, as $I_W$.

{\bf Step 2. General unbounded case.}
For any $L>0$, let $W^{(k)}_L:=W^{(k)}\wedge L$. Then $$H_{W_L}(\nu)=H(\nu|\alpha)+\sum_{k=2}^N \int_{S^k}W^{(k)}_L d\nu^{\otimes k},\ \text{ if }H(\nu|\alpha)<+\infty \text{ and $+\infty$ otherwise}$$
is inf-compact on $\MM_1(S)$, by Step 1. Therefore $H_W$ is also inf-compact since $H_{W_L}(\nu)\uparrow H_W(\nu)$ as $L \uparrow +\infty$.

For any closed subset $C$ of $\mathcal{M}_1(S),$ we have for any $L>0$,
\beq\label{ULDP1b}\aligned P_n^{*}\{L_n\in C\}&= \int_{S^n}1_{\{L_n\in C\}} \exp\left(-n\sum_{k=2}^N U_n(W^{(k)})\right)d\alpha^{\otimes n}\\
&\le \int_{S^n}1_{\{L_n\in C\}} \exp\left(-n\sum_{k=2}^N U_n(W^{(k)}_L)\right)d\alpha^{\otimes n} \\
&\le \exp\left(-n\inf_{\nu\in C}H_{W_L}(\nu)+o(n)\right),\\
\endaligned
\nneq
where the last inequality follows from Lemma \ref{LDP1} and Varadhan's Laplace lemma. Hence we get
$$\limsup_{n\to+\infty}\frac1n\log P_n^{*}\{L_n\in C\}\le -\inf_{\nu\in C}H_{W_L}(\nu),
$$
which leads to
\beq\label{ULDP2a}\limsup_{n\to+\infty}\frac1n\log P_n^*\{L_n\in C\}\le-\inf_{\nu\in C}H_W(\nu)
\nneq
because  $\inf_{\nu\in C}H_{W_L}(\nu) \uparrow \inf_{\nu\in C}H_{W}(\nu)$ as $L\uparrow+\infty$ by the inf-compactness of $H_{W_L}, H_W$.

Taking $C=\MM_1(S)$ in (\ref{ULDP2a}) we obtain
$$\limsup_{n\to+\infty}\frac1n\log\widetilde{Z}_n\le -\inf_{\nu\in\mathcal{M}_1(S)}H_W(\nu).
$$
 By the lower bound (\ref{LZ}) in Proposition  \ref{LLDP1} and the upper bound above,  and the fact that $H_W(\nu)=+\infty$ once if $W^{(k)}\notin  L^1(\nu^{\otimes k})$ for some $2\le k\le N$, we obtain
\beq\label{main1a}
\lim_{n\to\infty}  \frac 1n\log \tilde Z_n =-\inf_{\nu\in \MM_1(S)} H_W(\nu)
\nneq
which is a finite quantity (i.e. in $\rr$) by the inf-compactness of $H_W$ and (\ref{HWa}).
With  this key equality (\ref{main1a}) in hand,  the LDP of $\{P_n(L_n\in\cdot)\}_{n\ge N}$ follows from the lower bound in  Proposition  \ref{LLDP1} and the upper bound (\ref{ULDP2a}).

{\bf Step 3.}
Finally it remains to show the identification (\ref{rfc}) of the rate function $I_W$. For any $\nu$ such that $\nu\ll \alpha$ and $W^{(k),-}\in L^1(\nu^{\otimes k}),\ 2\le k\le N$, we have
$$
\aligned
\frac 1n H(\nu^{\otimes n}| P_n) &= \frac 1n \ee^{\nu^{\otimes n} } \left(\log \frac{d \nu^{\otimes n}} {d\alpha^{\otimes n} } + n\sum_{k=2}^N U_n(W^{(k)}) + \log \tilde Z_n\right)\\
&= H(\nu|\alpha) + \sum_{k=2}^N \int W^{(k)} d\nu^{\otimes n} + \frac 1n\log\tilde Z_n
\endaligned
$$
which yields (\ref{rfc}) by (\ref{main1a}), as $n\to +\infty$.

\nprf

\subsection{Proof of Theorem \ref{mainW}}
We first present the result of Sanov's theorem in the Wasserstein distance by Wang et. al. \cite{[WWW]}.
\bprop\label{SW} Let $\{X_n\}_{n\ge1}$ be a sequence of i.i.d. random variables defined on a probability space $(\Omega,\FF, P)$ with values in a Polish space $(S,\rho)$, of common law $\alpha$, then $\{P(L_n\in\cdot)\}_{n\ge1}$ satisfies the LDP on $(\MM_1^p(S),W_p)$ with speed $n$ and the good rate function $H(\cdot|\alpha )$, if and only if
$$\int_{S}\exp\{\lambda \rho^p(x,x_0)\}\alpha(dx)<+\infty,\ \ \forall \lambda>0$$
for some (hence for any) $x_0\in S$.
\nprop

\bprf[Proof of Theorem \ref{mainW}] Since we have established the LDP of the empirical measure $L_n$ under $P_n$ on $\MM_1(S)$ equipped with the weak convergence topology, it is sufficient to prove the exponential tightness of $\{P_n(L_n\in\cdot)\}_{n\ge N}$ in $(\MM_1^p(S),W_p)$ (see \cite[Corollary 4.2.6]{[DZ]}).

For any compact subset $K\subset \MM_1^p(S)$, and any fixed $a,b\in(1,+\infty)$ with $\frac1a+\frac1b=1$, we have by H\"older's inequality
\beq\label{ET3}\aligned P_n\{L_n\notin K\}&=\frac1{\widetilde{Z}_n}\int1_{\{L_n\notin K\}}\exp\left(-n\sum_{k=2}^N U_n(W^{(k)})\right)d\alpha^{\otimes n}\\
&\le \frac1{\widetilde{Z}_n}\left[\alpha^{\otimes n}\{L_n\notin K\}\right]^{1/a}\times\left(\int\exp\left(-bn\sum_{k=2}^N U_n(W^{(k)})\right)d\alpha^{\otimes n}\right)^{1/b}.\\
\endaligned
\nneq
Hence we have
\beq\label{ET2}\aligned &\limsup_{n\to+\infty}\frac1n\log P_n\{L_n\notin K\}\\
&\le \frac1a\limsup_{n\to+\infty}\frac1{n}\log\alpha^{\otimes n}\{L_n\notin K\}-\limsup_{n\to+\infty}\frac1n\log\widetilde{Z}_n\\
&\ \ \ \ \ \ \ \ \ \  \ \ \ \ \ \ \  \ \ \ \ \ \ \ \ \ \ \ \ \ \ +\frac1b\limsup_{n\to+\infty}\frac1{n}\log\int\exp\left(-n\sum_{k=2}^N U_n(bW^{(k)})\right)d\alpha^{\otimes n}\\
&= \frac1a\limsup_{n\to+\infty}\frac1{n}\log\alpha^{\otimes n}\{L_n\notin K\}+\inf_{\nu\in \MM_1(S)} H_W(\nu)-\frac1b\inf_{\nu\in \MM_1(S)} H_{bW}(\nu),\\
\endaligned
\nneq
where $\inf_{\nu\in \MM_1(S)} H_W(\nu)$ and $$\inf_{\nu\in \MM_1(S)} H_{bW}(\nu):=\inf_{\nu\in \MM_1(S)} \{H(\nu|\alpha)+\sum_{k=2}^N \int_{S^k}bW^{(k)}d\nu^{\otimes k}\}$$ are finite by (\ref{main1b}) in Theorem \ref{main}.

Note that under the exponential integrability condition (\ref{mainW1}), the LDP holds for $L_n$ under $\alpha^{\otimes n}$ with respect to the Wasserstein topology by Proposition \ref{SW}. Thus for any $L>0$, there exists a compact subset $K_L\subset \MM_1^p(S)$ such that
\beq\label{ET1} \limsup_{n\to+\infty}\frac1n\log\alpha^{\otimes n}\{L_n\notin K_L\}\le -aL-a\inf_{\nu\in \MM_1(S)} H_W(\nu)
+\frac{a}{b}\inf_{\nu\in \MM_1(S)} H_{bW}(\nu).
\nneq

Plugging (\ref{ET1})  into (\ref{ET2}), we get
\beq\label{ET4} \limsup_{n\to+\infty}\frac1n\log P_n\{L_n\notin K_L\}\le -L,
\nneq
which completes the proof.

\nprf

{\bf Acknowledgements}

The authors are grateful to the referees for their careful reading of the paper and for many valuable comments and suggestions.

\bibliographystyle{plain}

\end{document}